%% file: Fu5.tex
\documentclass[12pt]{article}
\usepackage{amsmath,amssymb,amsthm,amscd,amsfonts}
\usepackage[dvips]{graphicx}
\usepackage[dvips]{color}
\chardef\bslash=`\\
\newtheorem{Thm}{Theorem}[section]
\newtheorem{Cor}[Thm]{Corollary}
\newtheorem{Lem}[Thm]{Lemma}
\newtheorem{Prop}[Thm]{Proposition}

\newtheorem{Rque}[Thm]{Remark}
\newtheorem{Conj}{Conjecture}
\newtheorem{Pb}{Problem}

\title{Mukai flops and deformations of symplectic resolutions}
\author{Baohua Fu}

\def\cit{{\mathbb C}}
\def\nit{{\mathbb N}}

\def\pit{{\mathbb P}}

\def\0{{\mathcal O}}

\def\X{{\mathcal X}}
\def\Z{{\mathcal Z}}

\def\Y{{\mathcal Y}}
\def\W{{\mathcal W}}

\def\Sp{\mathop{\rm Sp}\nolimits}
\def\codim{\mathop{\rm codim}\nolimits}
\def\Spec{\mathop{\rm Spec}\nolimits}
\def\SL{\mathop{\rm SL}\nolimits}
\def\GL{\mathop{\rm GL}\nolimits}
\def\Hilb{\mathop{\rm Hilb}\nolimits}
\def\id{\mathop{\rm id}\nolimits}
\def\Fix{\mathop{\rm Fix}\nolimits}

\begin{document}
\maketitle
\begin{abstract}
We prove that two projective symplectic resolutions of $\cit^{2n}/G$ are connected by Mukai flops in codimension
$2$ for a finite sub-group $G < \Sp(2n)$. It is also shown that two projective symplectic resolutions of $\cit^4/G$
are deformation equivalent.
\end{abstract}
\section{Introduction}
A {\em symplectic variety} is a complex normal variety $W$ with a
holomorphic symplectic form on its smooth part which can be extended
to a global holomorphic form on any resolution. A resolution $Z \to
W$ of $W$ is called {\em symplectic} if the lifted holomorphic
$2$-form is again symplectic on $Z$.

 Examples of symplectic varieties
include the normalization of a nilpotent orbit closure in a
semi-simple complex Lie algebra and the quotient of $\cit^{2n}$ by a
finite subgroup $G < \Sp(2n)$.   The purpose of this paper is to
study projective symplectic resolutions of $\cit^{2n}/G$.

One way of constructing a symplectic resolution from another is to perform Mukai flops. This process can be
described as follows: let $W$ be a symplectic variety and $\pi: Z \to W$ a symplectic resolution. Assume that
$W$ contains a smooth closed subvariety $Y$ and that $\pi^{-1}(Y)$ contains a smooth subvariety $P$   such that the
restriction of $\pi$ to $P$ makes $P$ into a $\pit^l$-bundle over $Y$. If $\codim(P) = l$, we can blow up $Z$
along $P$ and then blow down along the other direction, which gives another (proper) symplectic resolution
$\pi^+: Z^+ \to W$. Notice that the resulting variety $Z^+$ may be not algebraic. Sometimes one needs to perform
simultaneously several Mukai flops to obtain a projective morphism $\pi^+$.
 The diagram $Z \to W
\leftarrow Z^+$ is called a {\em Mukai flop over $W$ with center $P$}. A {\em Mukai flop in codimension $2$} is a
diagram which becomes a Mukai flop after removing subvarieties of codimension greater than $2$.

As to the birational geometry in codimension $2$ of projective
symplectic resolutions, one has the following conjecture (due to
Hu-Yau \cite{HY}):
\begin{Conj}[Hu-Yau]
Any two projective symplectic resolutions of a symplectic variety
are connected by Mukai flops in codimension $2$.
\end{Conj}

This conjecture is true for four-dimensional symplectic varieties by the work of Wierzba and  Wi\'sniewski
(\cite{WW}) (partial results had previously  been obtained in \cite{BHL}, see also \cite{CMSB}) for the existence of
flops and by the work of Matsuki \cite{Mat} for the termination of flops. In \cite{Fu}, we have verified this
conjecture for symplectic resolutions of nilpotent orbit closures. The first result of this note is the
following theorem (partial but stronger results had been proven in \cite{Fu}).
\begin{Thm}
Let $G < \Sp(2n)$ be a finite subgroup.
Any two projective symplectic resolutions of $\cit^{2n}/G$ are
connected by Mukai flops in codimension $2$.
\end{Thm}

The  idea is to reduce the problem to dimension $4$, and then apply \cite{WW}.
The main technics (see Proposition \ref{product} and Lemma \ref{local})  for this reduction are already contained in \cite{Ka1}.
The proofs here are slightly different.

 Then we
study deformations of symplectic resolutions. Recall that a {\em
deformation} of a variety $X$ (usually not compact) is a flat
morphism $\X \xrightarrow{q} S$ from a variety $\X$ to a pointed
smooth connected curve $0 \in S$ such that $q^{-1}(0) \simeq X$. A
deformation of a proper morphism $X \xrightarrow{f} Y$ is an
$S$-morphism $\X \xrightarrow{F} \Y$ such that $\X \to S$ (resp. $\Y
\to S$) is a deformation of $X$ (resp. $Y$) and $F_0 = f.$

Let $ X \xrightarrow{f} Y \xleftarrow{f^+} X^+$ be two proper
morphisms. One says that $f$ and $f^+$ are {\em deformation
equivalent} if there exist deformations of $f$ and $f^+$: $\X
\xrightarrow{F} \Y \xleftarrow{F^+} \X^+$ such that for a general
point $s \in S$ the morphisms $\X_s \xrightarrow{F_s} \Y_s
\xleftarrow{F^+_s} \X^+_s$ are isomorphisms. Motivated by
 results of D. Huybrechts (\cite{Huy2}), it is conjectured
in \cite{FN} (see also \cite{Ka2}) that
\begin{Conj}\label{deform}
 Any two symplectic resolutions
of a symplectic variety are deformation equivalent.
\end{Conj}

For nilpotent orbit closures of classical type, this conjecture is
proved by Y. Namikawa in \cite{Nam} (the case of
$\mathfrak{sl}(n)$ had previously been  proved in \cite{FN}). In \cite{Ka2}, D.
Kaledin constructed the so-called {\em twister deformation} of a
symplectic resolution (under mild assumptions).
Combining this with a trick of D. Huybrechts (\cite{Huy}) and results
in \cite{WW}, we prove the following
\begin{Thm}
Let $G < \Sp(4)$ be a finite subgroup.
Any two projective symplectic resolutions of $\cit^4/G$ are
deformation equivalent.
\end{Thm}

This note ends with a study of symplectic resolutions of the wreath
product $W: = (\cit^2/\Gamma)^{(n)}$, where $\Gamma < \SL(2)$ is a
finite subgroup. It is conjectured that any two projective symplectic
resolutions of $W$ are connected by Mukai flops with flop center
contained in the fiber over $0 \in W$. In the case of $\Gamma$ being
of type $A_k$ and $n=2$, we give a way to describe  all possible
projective symplectic resolutions of $W$. \\

{\em Acknowledgments}: I am glad to thank  S. Druel for the
remarks, corrections and suggestions that he   made to a first
version of this article, and  O. Debarre for a careful reading.

\section{Birational geometry in codimension $2$}
We begin with the following proposition, which is proved (as is Lemma 2.4 later) in the formal setting by D. Kaledin
(\cite{Ka1} Proposition 5.2).
\begin{Prop} \label{product}
Let $W$ be a symplectic variety and $\Delta^{2l}$ the open unit disk in $\cit^{2l}$.
Then any projective symplectic resolution of
$W \times \Delta^{2l}$ is of the form $Z \times \Delta^{2l}
\xrightarrow{\pi} W \times \Delta^{2l}$, where $Z \xrightarrow{\pi'}
W$ is a symplectic resolution and $\pi = \pi' \times \id$.
\end{Prop}
\begin{proof}
Suppose that we have a  symplectic resolution $X \xrightarrow{\pi}
W \times \Delta^{2l}$. For any non-zero vector  $v \in \Delta^{2l}$,
it defines a constant vector field $\mathfrak{t}_v$ on the smooth
part, say $U$ of $W \times \Delta^{2l}$. Furthermore, on $U$, one
has an isomorphism of sheaves $\Omega^1 \simeq \mathcal{T}$, under
which the vector field $\mathfrak{t}_v$ corresponds to a 1-form
$\alpha_v$. It is easy to show that $\alpha_v = p_2^* \beta$ for
some 1-form $\beta$ on $\Delta^{2l}$, where $p_2: W \times
\Delta^{2l}\to \Delta^{2l}$ is the projection to the second factor. In
particular,  $\alpha_v$ extends to a well-defined 1-form on the
whole of  $W \times \Delta^{2l}$. Let $\mathfrak{t}'_v$ be the
vector field on $X$ corresponding to the 1-form $\pi^* \alpha_v$
under  the isomorphism  $\Omega^1_X \simeq \mathcal{T}_X$. Then
$\mathfrak{t}'_v$ is the vector field lifting $\mathfrak{t}_v$.
Furthermore $\mathfrak{t}'_v$ vanishes nowhere on $X$, thus it
defines a holomorphic flow $\phi_v(t)$ on $X$ (see the proof of Theorem 1.3 \cite{Ka1}).

Let $q: X \to \Delta^{2l}$ be the composition $p_2 \circ \pi$ and $Z
= q^{-1}(0). $ Let $\pi': Z \to W$ be the restriction of $\pi$ to
$Z \to W \times \{0 \}$. Then the flow $\phi_v(t)$ satisfies
$q(\phi_v(t)(z)) = t v$ for any $z \in Z$. We define a morphism $Z
\times \Delta^{2l} \to X$ as follows: $(z, v) \mapsto \phi_v(1)(z).$
One sees easily that this is an isomorphism. Moreover one has
$\pi(\phi_v(1)(z)) = (\pi'(z), v)$.

In conclusion, we obtain a decomposition $X = Z \times \Delta^{2l}$, a
map $\pi': Z \to W$ and an isomorphism $\pi = \pi' \times \id$. That
$\pi$ is a symplectic resolution implies
 that $Z$ is smooth and  $\pi'$ is a  symplectic resolution of $W$.
\end{proof}

The same arguments hold if one replaces $\Delta^{2l}$ by $\cit^{2l}$.   An immediate corollary is the
following (which is also proved in \cite{Ka1} Theorem 1.6):
\begin{Cor}
Let $V_i$ be a symplectic vector space and $G_i < \Sp(V_i)$ a finite
subgroup, $i \in \{1, 2\}$. Then $V_1/G_1 \times V_2/G_2$ admits  a
 symplectic resolution if and only if $V_1/G_1$ and $
V_2/G_2$ both admit  symplectic resolutions.
\end{Cor}
\begin{proof}
Take a smooth point $v \in V_1/G_1$ and a neighborhood isomorphic to the unit disk $\Delta$. If the product
admits a symplectic resolution, so does $\Delta \times V_2/G_2$. The precedent proposition then implies that
$V_2/G_2$ admits a symplectic resolution. Similarly $V_1/G_1$ also admits a symplectic resolution.
\end{proof}

\begin{Rque} \upshape
We do not know if every projective symplectic resolution of $V_1/G_1
\times V_2/G_2$ is a product of resolutions of $V_1/G_1$ and $
V_2/G_2$. This is true if $G_1$ or $G_2$ is trivial by the precedent
proposition.
\end{Rque}

From now on, let $V$ be a $2n$-dimensional symplectic vector space
and $G < \Sp(V)$ a finite subgroup. We denote by $W$ the quotient
space $V/G$. We have the  rank stratification on $W$ defined as $V_k
= \{v \in V | \codim  V^{G_v} = 2k \}$. The quotient $W_k = V_k /G$ is a
smooth algebraic variety of dimension $2n - 2k$ and $W_0$ is the
smooth part of $W$. Moreover, the projection $V_k \to W_k$ is
\'etale (Lemma 4.1 \cite{Ka1}).

Take a component $Y $ of $W_k$ and a connected component  $V_Y$ of
the preimage of $Y$ in $V_k$. Let $H$ be the stabilizer of a point
in $V_Y$.  Then $H$ is independent of the choice of the point and
$V_Y$ is a Zariski open set in the $H$-fixed subspace  $V^H$. Let
$N(H)$ be the normalizer of $H$ in $G$ and $Q(H) = N(H)/H$ the
quotient group. One shows that $V^H$ is $N(H)$-invariant. Since $H$
acts trivially on $V^H$, one obtains an action of $Q(H)$ on $V^H$,
which is a free action on $V_Y$ and we have an isomorphism $V_Y/Q(H)
\simeq Y.$

Let $V_H$ be the annihilator of $V^H$ with respect  to the symplectic form $\omega_0$  on $V$. Notice that
$\omega_0$ restricted to $V^H$ is again symplectic, thus one has a decomposition $V = V^H \oplus V_H$, which is
$N(H)$-invariant. Furthermore $N(H)$ acts symplectically on $V_H$, i.e., $N(H) < \Sp(V_H)$. This decomposition
induces a morphism $\mu: (V^H \times V_H/H)/Q(H) \to V/G$  which maps $V_Y/Q(H) \times \{0\} $ isomorphically to
$Y$ and  $\mu$ is \'etale in a Zariski open set containing $V_Y/Q(H) \times \{0\}$. For more details, see
section 4 of \cite{GK}.  This implies (see also Lemma 4.2 \cite{Ka1}):
\begin{Lem}\label{local}
Any point in $Y$ admits an analytical open neighborhood which is isomorphic to $\Delta^{2l} \times D_H$, where
$\Delta^{2l}$ is the unit disk of dimension $\dim (V^H)$ and $D_H$ is the image of the unit disk in $V_H$ under
the projection $V_H \to V_H/H$.
\end{Lem}

\begin{Rque} \label{dim4} \upshape
In the case of $\dim(V_H) = 4$,  any symplectic resolution of $D_H$ extends to a symplectic resolution of $S:=
V_H/H$. In fact, outside the zero point, $S$ has only $ADE$ singularities and
 $S - \{0\}$ admits a unique symplectic resolution $\tilde{S} \to S-\{0\}$.
 Now any symplectic resolution of $D_H$ agrees  automatically with $\tilde{S}$ over
 $D_H - \{0\}$, thus it pastes with $\tilde{S}$ to a symplectic
 resolution of $S$. If the resolution of $D_H$ is projective, then the one obtained for
$V_H/H$ is again projective.
\end{Rque}

Suppose that we have two projective symplectic resolutions  $Z
\xrightarrow{\pi} W  \xleftarrow{\pi^+} Z^+$. Let $\phi$ be the
rational map $\pi^{-1} \circ \pi^+: Z^+ \dashrightarrow Z$.
\begin{Lem} \label{isocodim1}
The rational map $\phi$ induces an isomorphism from $(\pi^+)^{-1}(U)$ to $\pi^{-1}(U)$,
where $U = W_0 \cup W_1$.
\end{Lem}
\begin{proof}
By the   lemma above, every point $y \in W_1$ admits a
neighborhood $U_y$ isomorphic to $\Delta^{2n-2} \times D_H$ for
some finite subgroup $H < \SL(2)$. By Proposition \ref{product},
every symplectic resolution of $U_y$ is a product of
$\Delta^{2n-2}$ with a symplectic resolution of $D_H$, while $D_H$
admits a unique symplectic resolution given by the minimal
resolution, thus $\phi$ is an isomorphism  from
$(\pi^+)^{-1}(U_y)$ to $\pi^{-1}(U_y)$.
\end{proof}
\begin{Thm}
Two projective symplectic resolutions of $W$ are connected by Mukai flops over $W$ in codimension $2$.
\end{Thm}
\begin{proof}
Let $Z \xrightarrow{\pi} W  \xleftarrow{\pi^+} Z^+$ be two
projective symplectic resolutions. By the semi-smallness of
symplectic resolutions (Prop. 4.4 \cite{Ka1}),
$\pi^{-1}(\overline{W_3})$ (respectively
$(\pi^+)^{-1}(\overline{W_3})$) has codimension at least $3$ in $Z$
(resp. $Z^+$). Since we are  interested in the codimension $2$
birational geometry, we can replace $W$ by $W_0 \cup W_1 \cup W_2$.
By the precedent lemma, $\phi$ is already an isomorphism over $W_0
\cup W_1$.

Take a connected component $Y$ in $W_2$ and a point $y \in Y$.
Then there exists an analytical neighborhood $U_y$ of $y$
isomorphic to $\Delta^{2n-4} \times D_H$ for some finite subgroup
$H < \Sp(4)$. By proposition \ref{product}, the projective
symplectic resolution $\pi^{-1}(U_y) \to U_y$ is isomorphic to the
product $\Delta^{2n-4} \times X \to U_y$, where $X \to D_H$ is a
projective symplectic resolution. Similarly for $\pi^+$, one finds
another projective symplectic resolution $X^+ \to D_H$. Since $Y$
is connected, $X$,  $X^+$ and their morphisms to $D_H$ are
independent of the choice of $y$. By Remark \ref{dim4}, these two
symplectic resolutions come in fact from symplectic resolutions of
$\cit^4 / H$.

By \cite{WW} and \cite{Mat}, the birational map $X \dashrightarrow X^+$ is decomposed as a sequence of Mukai flops.
Without any loss of generality, one may suppose that $X  \dashrightarrow  X^+$ is a Mukai flop with flop center $P
\subset X$. Since $X$ is independent of the choice of $y$, one can find a subvariety $E$ in $Z$ which is a
fibration over $Y$ with fibers isomorphic to  $P$.
By the McKay correspondence (see \cite{Ka3}), irreducible components of
codimension $2$ in $\pi^{-1} (Y)$ correspond to dimension $2$ components in the central fiber of $X$.
The subvariety $E$ is then the irreducible component of codimension $2$ in the preimage of $Y$ corresponding to $P$.

Now if we perform a Mukai flop in $Z$ along $E$, one obtains
another symplectic resolution $X' \xrightarrow{\pi'} W$ such that the rational map $X'  \dashrightarrow  X^+$ is an
isomorphism between preimages of  $Y$.

Now if  we do the same operations for other components in $W_2$,
one arrives finally at the resolution $\pi^+$.
\end{proof}

We end this section by the following proposition, whose proof is clear.
\begin{Prop}
Let $W:= V/G$ be a quotient symplectic variety. Suppose that for
every component $Y$ in $W_2$, the corresponding $4$-dimensional
quotient $\cit^4/H_Y$ admits a unique projective symplectic
resolution. Then any two projective symplectic resolutions of $W$
are isomorphic in codimension $2$.
\end{Prop}

The following $\cit^4/G$   admit a unique projective symplectic
resolution:

(i) $\cit^2/G_1 \times \cit^2/G_2$ where $G_1,
G_2$ are finite subgroups of $\SL(2)$;

(ii) $(T^*\cit^2)/G$, where $G
< \GL(2)$ such that $\{g | \Fix(g) = 0 \}$ form a single conjugacy
class.

Case (i) follows from \cite{WW} since the central fiber
contains no copies of $\pit^2$, while case (ii) is proved in
\cite{Fu} (Cor. 1.3).

\section{Deformation equivalence}

Let $V$ be a $2n$-dimensional symplectic vector space and $G <
\Sp(V)$ a finite subgroup. Suppose that we have a projective
symplectic resolution $\pi: Z \to W: = V/G$. Take a $\pi$-ample line
bundle $L$ on $Z$. By \cite{Ka2}, there exists a twister deformation
of $\pi$ over the formal disk $\Spec(\cit[[x]])$. Since $W$ admits an
expanding $\cit^*$-action (i.e., positively weighted)
which lifts to $Z$ via $\pi$, this twister
deformation extends to an actual deformation over $S = \cit$, say
$\mathcal{Z} \xrightarrow{\Phi} \mathcal{W}$ (see Lemma A. 15 and
Proposition 5.4 \cite{GK}). Furthermore, for a generic $s \in S$,
the morphism $\Phi_s: \mathcal{Z}_z \to \mathcal{W}_s$ is an
isomorphism. Moreover, by \cite{Ka2}, the Kodaira-Spencer class $v$
of the deformation $\mathcal{Z} \to S$ is nothing but $c_1(L) \in
H^1(Z, T_Z) \simeq H^1(Z, \Omega_Z).$

Let $P \subset Z$ be a subvariety isomorphic to $\pit^n$. Denote by
 $\bar{v}$ the image  of the Kodaira-Spencer class $v$ under the morphism
$H^1(Z, \Omega_Z) \to H^1(P, \Omega_P)$. The following lemma is a
special case of Lemma 3.6 \cite{Huy}. We omit the proof here.
\begin{Lem}\label{normal}
If $\bar{v}$ is non-zero, then $\mathcal{N}_{P|\mathcal{Z}} \simeq
\0_P(-1)^{\oplus n+1}$.
\end{Lem}

Let $p: \tilde{\Z} \to \Z$ be the blow up of  $\mathcal{Z}$ along $P$. Under the assumption of the precedent
lemma, the exceptional divisor $E$ is isomorphic to $\pit(\0_P(-1)^{\oplus n+1}) = P \times P^*$, where $P^*$ is
the dual of $P$, and the normal bundle $\mathcal{N}_{E|\Z} $ is the tautological bundle.  In particular, the
restriction of $\0_{\tilde{\Z}}(E))$ to any fiber of $P \times P^* \to P^*$ is $\0(-1)$. By Nakano-Fujiki
criterion, there exists a contraction $\tilde{\Z} \to \Z^+$ which blows down $E$ to $P^*$. Let $Z^+$ be the
Mukai flop of $Z$ along $P$. Then $\Z^+$ is a one-parameter deformation of $Z^+$.
 Let $L^+$ be the
strict transform of $L$ under the rational map $Z  \dashrightarrow Z^+$.
\begin{Lem}\label{KSclass}
$c_1(L^+)$ is the Kodaira-Spencer class of the deformation $\Z^+ \to
S$.
\end{Lem}
\begin{proof}
Let $U = Z - P$, isomorphic to  $U^+: = Z^+ - P^*$. We denote by $v|_U$ (resp.
$v^+|_{U^+}$) the image of the Kodaira-Spencer class under the map
$H^1(Z, T_Z) \to H^1(U, T_U)$ (resp. $H^1(Z^+, T_{Z^+}) \to H^1(U^+,
T_{U^+})$). Notice that we have an $S$-isomorphism $\Z - P \simeq
\Z^+ -P^*$; thus $v|_U = v^+|_{U^+} $ via the isomorphism $U \simeq
U^+$.

The map $H^1(Z, T_Z) \to H^1(U, T_U)$ is injective since $\codim_Z P
\geq 2$. Furthermore $v|_U = c_1(L)|_U = c_1(L^+)|_{U^+}$, thus
$c_1(L^+)$ is the Kodaira-Spencer class of the deformation $\Z^+ \to
S$.
\end{proof}

If furthermore $P$ is mapped to a point by $\pi$, then  one has
another symplectic resolution $Z^+ \to W$  which admits a
deformation $\Z^+ \to \W$. The deformations  one wants to construct
in the following theorem are based on this.

\begin{Thm}
Let $V$ be a four-dimensional symplectic vector space and $G <
\Sp(V)$ a finite subgroup. Then any two projective symplectic
resolutions of $V/G$ are deformation equivalent.
\end{Thm}
\begin{proof}
Let $W = V/G$ and $Z \xrightarrow{\pi} W \xleftarrow{\pi^+} Z^+$ two projective symplectic resolutions.  Take a
$\pi^+$-ample line bundle $L^+$ on $Z^+$.
 Then we have a deformation of $\pi^+$: $\Z^+ \to \W$ such that $c_1(L^+)$
is the Kodaira-Spencer class of $\Z^+ \to S$. Let $L$ be the strict transform to $Z$ of $L^+$. Then $L$ is
$\pi$-big. If $L$ is $\pi$-nef, then the two resolutions $\pi$ and $ \pi^+$ are isomorphic (see \cite{FN} Theorem 2.2).

If $L$ is not $\pi$-nef, we can find a $(Z, \epsilon L)$-extremal ray $R$ for small $\epsilon > 0$ (see \cite{KMM}).
The locus $E$ of $R$ in $Z$ is contained in $\pi^{-1}(0)$ by Lemma \ref{isocodim1} and the contraction of $R$
gives a small contraction since $dim (\pi^{-1}(0)) \leq 2$ by the semi-smallness of symplectic resolutions. By
\cite{WW}, $E$ is a disjoint union of copies isomorphic to $\pit^2$. Furthermore $L$ is negative on every curve
in $E$. We can perform a Mukai flop along $E$ to obtain $\pi_1: Z_1 \to W$. The strict transform $L_1$ of $L$ is
then positive on all curves of $E^*$.  If $L_1$ is not $\pi_1$-nef, then we can continue this process. After
finitely many steps, say $Z  \dashrightarrow Z_1  \dashrightarrow \dots \to Z_{l+1}$ one arrives to $\pi_{l+1} = \pi^+$.

Let $L_i$ be the strict transform of $L$ to $Z_i$ and $E_i$ the flop
center of $Z_i  \dashrightarrow  Z_{i+1}$. Then $L_{i+1}$ is positive on curves
in $E_i^*$ for $i=1, \cdots, l$.  By Lemma \ref{normal}, the normal
bundle $N_{E_l^*|\Z^+}$ is isomorphic to $\0_{E_l^*}(-1)^{\oplus
3}$. Thus we can blow up $\Z^+$ at $E_l^*$ then blow down along the
other direction to obtain a deformation of $\pi_l$: $\Z_l \to \W$.
By Lemma \ref{KSclass} and Lemma \ref{normal}, one can perform the
same process to $E_{l-1}^*$ in $\Z_l$ and so on. Finally one obtains
a deformation of $\pi$: $\Z \to \W$. Then the two deformations $\Z
\to \W \leftarrow \Z^+$ give the equivalence.
\end{proof}

\section{Wreath product and Hilbert schemes}

Let $\Gamma < \SL(2)$ be a finite subgroup and $W =
(\cit^2/\Gamma)^{(n)}$ the $n$-th symmetric product of
$\cit^2/\Gamma$. Then $W$ is the quotient of $\cit^{2n}$ by the
wreath product $\Gamma_n = \Gamma \sim \mathcal{S}_n$. Explicitly,
$\Gamma_n = \{(g, \sigma) | g \in \Gamma^n, \sigma \in \mathcal{S}_n
\}$ with the multiplication $(g, \sigma)\cdot  (h, \tau) = (g
\sigma(h), \sigma \tau)$, where $\sigma(h) = (h_{\sigma^{-1}(1)},
\ldots, h_{\sigma^{-1}(n)})$.

Let $S \to \cit^2/\Gamma$ be the minimal resolution. Then the
composition $$\Hilb^n(S) \xrightarrow{\tau} S^{(n)} \to
(\cit^2/\Gamma)^{(n)}$$ gives a projective symplectic resolution
$\Hilb^n(S) \xrightarrow{\pi} W$ (see also \cite{Wan}). When $\Gamma$
is trivial, this is the unique projective symplectic resolution of
$W$ (cf. \cite{FN}). However it is not true for a non-trivial
$\Gamma$. The following problem is open for $n \geq 3$.
\begin{Pb}
Find out all projective (resp. proper) symplectic resolutions of $W
= \cit^{2n}/\Gamma_n$.
\end{Pb}

Let $C_i, i \in \{1, \ldots, k\}$ be the irreducible components in the
exceptional divisor $S \to \cit^2/\Gamma$. Then in the central fiber
$\pi^{-1}(0)$ there are $k$ disjoint copies of $\pit^n$, given by
the strict transforms of $C_i^{(n)}$ via $\tau$. In particular, we
can perform Mukai flops to obtain many different symplectic
resolutions of $W$. However it is not clear if these resolutions are
still projective. An answer to Problem 1 is expected from
the following
\begin{Conj}
Any two projective symplectic resolutions of $W$ are connected by a
sequence of Mukai flops with flop centers contained in the fiber over
$0 \in W$. In particular, they are isomorphic over $W - 0$.
\end{Conj}

It is not totally unlikely that the precedent conjecture holds for
any quotient variety $V/G$ which is not a product of quotient
varieties. For $4$-dimensional quotients, this is true thanks to
the results in \cite{WW}.

A positive answer to this conjecture may imply that
Conjecture \ref{deform} is valid for projective symplectic
resolutions of $W$, by the arguments of the precedent section and
results in \cite{CMSB}.

There exists another natural symplectic resolution of $W$ that is constructed as follows
(constructed in \cite{Wan}): let $N = |\Gamma|$ be the order of
$\Gamma$. The action of $\Gamma$ on $\cit^2$ extends to
$\Hilb^{nN}(\cit^2)$ and $(\cit^2)^{(nN)}$. Thus the Hilbert-Chow
morphism induces a morphism between $\Gamma$-fixed points $\Hilb^{nN,
\Gamma}(\cit^2) \to (\cit^2)^{(nN), \Gamma}$. Notice that
$(\cit^2)^{(nN), \Gamma}$ is naturally identified with $W =
(\cit^2/\Gamma)^{(n)}$. Let $Z_{\Gamma, n}$ be the closure in
$\Hilb^{nN, \Gamma}(\cit^2)$ of unordered $n$-tuple of distinct
$\Gamma$-orbits in $\cit^2 - 0$. It is shown in \cite{Wan} that
$Z_{\Gamma, n}$ is a connected component of $\Hilb^{nN,
\Gamma}(\cit^2)$, thus it is smooth and symplectic. Moreover, the
morphism $ Z_{\Gamma, n} \xrightarrow{\pi^+}   (\cit^2)^{(nN),
\Gamma} \simeq W$ is an isomorphism over $W_{reg}$, thus it gives a
projective symplectic resolution of $W$.
\begin{Pb}
Connect the two resolutions  $\pi, \pi^+$ by Mukai flops.
\end{Pb}
\begin{Rque} \upshape
$\pi$ and $ \pi^+$ are in general non-isomorphic. In the case of
$\Gamma = \{\pm 1\}$ and $n= 2$, $\pi$ and $ \pi^+$ are the two
non-isomorphic projective symplectic resolutions that $(\cit^2/{\pm
1})^{(2)}$ can have (see \cite{FN}, Example 2.7).
\end{Rque}

In the following we give a way to describe all possible projective
symplectic resolutions of $W = (\cit^2/\Gamma)^{(2)}$. The
irreducible components in $\pi^{-1}(0)$ can be described as follows:

(i) $P_{i,i}$ ($1 \leq i \leq k$): the strict transform of
$C_i^{(2)}$ via $\tau$. They are isomorphic to $\pit^2$;

(ii) $P_{i,j}$ ($1 \leq i < j  \leq k$): the strict transform via
$\tau$ of the image of $C_i \times C_j$ under the morphism $S^2 \to
S^{(2)}$. If $C_i \cap C_j = \emptyset$, then $P_{i,j}$ is
isomorphic to $\pit^1 \times \pit^1$. If $C_i \cap C_j = \{x\}$,
then $P_{i,j}$ is isomorphic to the one point blow up of $\pit^1
\times \pit^1$;

(iii) $Q_i$ ($1 \leq i \leq k$): the preimage
$\tau^{-1}(\Delta_{C_i})$, where $\Delta_{C_i}$ is the diagonal
embedding  of $C_i$ in $S^{(2)}$. It is isomorphic to
$\pit(T_S|_{C_i}) \simeq \pit(\0_{C_i}(2) \oplus \0_{C_i}(-2))$,
thus a Hirzebruch surface $F_4$.
\begin{Lem}
The strict transform of $Q_i$ under any sequence of Mukai flops
along components in $\pi^{-1}(0)$ is not isomorphic to $\pit^2$.
\end{Lem}
\begin{proof}
To simplify the presentation, we will only prove the lemma for
$\Gamma$ being of type $A_k$, i.e., $\Gamma$ is a cyclic subgroup in
$\SL(2)$ of order $k+1$.  Let $C_i \cap C_{i+1} = \{x_i\}$ for $i =1,
\cdots, k$. One checks that $l_i: = Q_i \cap P_{i,i}$ is a conic in
$P_{i,i}$ and a negative section in $Q_i$. If we perform a Mukai
flop along $P_{i,i}$, then $l_i$ is transformed to a conic in
$P_{i,i}^*$, which is still called the strict transform of $l_i$.
The strict transform of $Q_i$ is isomorphic to $Q_i$. Among
$P_{i,j}$, only $P_{i-1,i}$ and $P_{i,i+1}$ intersect $l_i$, both
at one point (with multiplicity $2$).

One way to make the self-intersection of the strict transform of $l_i$
 positive is to flop  $P_{i-1,i}$ or $P_{i,i+1}$.
To do so, one needs to flop $P_{i,i}$ first. After the flop along $P_{i,i}$,
$P'_{i-1,i}$ intersects $P_{i,i}^*$ at one point (which lies on the strict transform of $l_i$).
By this, one sees that  the self-intersection of the strict
transform of $l_i$ is always negative.
\end{proof}

Thus to construct Mukai flops, one only needs to consider the
components $P_{i,j}$.  In the following we will assume that
$\Gamma$ is of type $A_k$ (with minor changes, analogue results can
be obtained for types $D_k, E_l$). The configuration of $P_{i,j}$
will be represented in $\nit^2$ as follows: $P_{i,j}$ is placed at
the position $(i,j)$, represented by a rectangle (resp. an ellipse,
a $\oplus$, a circle) if $P_{i,j}$ is isomorphic to $\pit^1 \times
\pit^1$ (resp. one point blow-up of $\pit^1 \times \pit^1$,
Hirzebruch surface $F_1$, $\pit^2$). These are the vertices of the
graph. It is easy to see that the intersection of components of
$P_{i,j}$ is either one point or a $\pit^1$ if not empty. Two
vertices are joined by a solid line (resp. dotted line) if their
intersection is a $\pit^1$ (resp. a point).

When we perform a Mukai flop at a vertex say $P_{i,j}$, the solid
line (resp. dotted line) joining this vertex is replaced by a dotted
line (resp. solid line). Other lines are untouched except the
following case: the vertex $P_{i,j}$ is joined to two vertices $P_1,
P_2$ by dotted lines. Then after the flop, the two dotted lines are
replaced by solid lines, and furthermore $P_1$ and $P_2$ are joined
by a dotted line.  Surely this process is symmetric, i.e., if
$P_{i,j}$ is joined to $P_1, P_2$ by solid lines and $P_1, P_2$ are
joined by dotted line, then after the flop along $P_{i,j}$, the
dotted line between $P_1$ and $P_2$ should be removed, and the solid
lines joining $P_{i,j}$ to $P_1, P_2$ are replaced by dotted ones.

Now we describe how the vertex labels change.  Since the process is
symmetric, we only describe the changes when  $P$ is a vertex
joined to $P_{i,j}$ by a solid line. Suppose that $P$ is labeled by
an ellipse (i.e., a one point blow-up of $\pit^1 \times \pit^1$). There
are two cases:

(i) the solid line comes from a dotted line, i.e., this line
corresponds to the exceptional fiber of the one point blow-up of
$\pit^1 \times \pit^1$, then the label at $P$ is changed to a square
(i.e., $\pit^1 \times \pit^1$) ;

(ii) otherwise, the label at $P$ is changed to be $\oplus$ (i.e.,
$F_1$).

If $P$ is labeled by a $\oplus$, then it is changed to a circle.
The following pictures are examples of symplectic resolutions of
$\cit^4/\Gamma_2$ with $\Gamma$  of type $A_4$.

Any projective symplectic resolutions of $W$ is obtained in this
way. However, it is not clear (and it may be not true) that any
sequence of Mukai flops gives a projective symplectic resolution.
Sometimes one needs to flop simultaneously  several disjoint
$\pit^2$  to obtain a projective resolution.

\newpage
$$
\begin{picture}(300,0)
\put(-60,-500){\input{model.pstex_t} }
\end{picture}
$$

 \vspace{0.4 cm}

 Labo. J. Leray,
 Facult\'e  des sciences, Universit\'e de Nantes

2, rue de la Houssini\`ere,  BP 92208 F-44322 Nantes Cedex 03 -
France

fu@math.univ-nantes.fr

\end{document}

%% file: model.pstex_t
\begin{picture}(0,0)%
\includegraphics{model.pstex}%
\end{picture}%
\setlength{\unitlength}{2565sp}%
\begingroup\makeatletter\ifx\SetFigFont\undefined%
\gdef\SetFigFont#1#2#3#4#5{%
  \reset@font\fontsize{#1}{#2pt}%
  \fontfamily{#3}\fontseries{#4}\fontshape{#5}%
  \selectfont}%
\fi\endgroup%
\begin{picture}(10561,12216)(289,-11815)
\end{picture}

%% file: Fu5.bbl
\newpage
\begin{thebibliography}{10}
\bibitem[BHL]{BHL}
Burns, D.; Hu, Y.; Luo, T., \emph{HyperK\"ahler manifolds and
birational transformations in dimension $4$}, Vector bundles and
representation theory (Columbia, MO, 2002),  141--149, Contemp.
Math., {\bf 322}, Amer. Math. Soc., Providence, RI, 2003.
\bibitem[CMSB]{CMSB}
Cho, K.; Miyaoka, Y.; Shepherd-Barron, N. I.,
\emph{Characterizations of projective space and applications to
complex symplectic manifolds}, Higher dimensional birational
geometry (Kyoto, 1997), 1--88, Adv. Stud. Pure Math., {\bf 35},
Math. Soc. Japan, Tokyo, 2002.
\bibitem[Fu]{Fu}
Fu, B., \emph{Birational geometry in codimension 2 of symplectic
resolutions},  math.AG/0409224.
\bibitem[FN]{FN}
Fu, B.; Namikawa, Y., \emph{Uniqueness of crepant resolutions and sympletic singularities},
Ann. Inst. Fourier, {\bf 54} (2004), no. 1, 1--19.
\bibitem[GK]{GK}
Ginzburg, V.;  Kaledin D., \emph{ Poisson deformations of symplectic quotient singularities},
Adv. Math. {\bf 186} (2004), no. 1, 1--57.
\bibitem[HY]{HY}
Hu, Y.; Yau, S.-T., \emph{HyperK\"ahler manifolds and birational
transformations},  Adv. Theor. Math. Phys.  {\bf 6}  (2002),  no. 3, 557--574.
\bibitem[Huy1]{Huy}
 Huybrechts, D., \emph{ Birational symplectic manifolds and their
 deformations},  J. Differential Geom. {\bf  45}  (1997),  no. 3, 488--513.
\bibitem[Huy2]{Huy2}
Huybrechts, D., \emph{ Compact hyper-K\"ahler manifolds: basic
results}, Invent. Math.  {\bf 135}  (1999),  no. 1, 63--113.
\bibitem[Ka1]{Ka1}
Kaledin, D., \emph{ On crepant resolutions of symplectic quotient
singularities},  Selecta Math. (N.S.)  {\bf 9}  (2003),  no. 4, 529--555.
\bibitem[Ka2]{Ka2}
Kaledin, D., \emph{Symplectic resolutions: deformations and
birational maps},  math.AG/0012008.
\bibitem[Ka3]{Ka3}
 Kaledin, D., \emph{McKay correspondence for symplectic quotient singularities},
Invent. Math.  {\bf 148}  (2002),  no. 1, 151--175.
\bibitem[KMM]{KMM}
Kawamata, Y., Matsuda, K., Matsuki, K., \emph{Introduction to the minimal model program}, in {\em Algebraic
Geometry, Sendai (1985)}, Advanced Studies in Pure Math. {\bf 10} (1987), Kinokuniya and North-Holland, 283-360.
\bibitem[Mat]{Mat}
Matsuki, K., \emph{Termination of flops for 4-folds}, Amer. J. Math. {\bf 113} (1991), 835-859.
\bibitem[Nam]{Nam}
Namikawa, Y., \emph{Birational Geometry of symplectic resolutions of
nilpotent orbits}, math.AG/0404072.
\bibitem[Wan]{Wan}
Wang, W., \emph{Hilbert schemes, wreath products, and the McKay
correspondence}, math.AG/9912104.
\bibitem[WW]{WW}
Wierzba, J.; Wi\'sniewski, J., \emph{Small contractions of symplectic
4-folds},
Duke Math. J. {\bf 120} (2003), no. 1, 65--95.
\end{thebibliography}
